\documentclass[12pt]{amsart}
\usepackage[english]{babel}

\usepackage{amssymb,amsmath,amscd,amsthm}

\def\{{\protect\lbrace}
\def\}{\protect\rbrace}

\begin{document}

\begin{center}
\textbf{\Large Right Serial Skew Laurent Series Rings}

A. A. Tuganbaev
\end{center}
\hfill National Research University "MPEI"

\hfill e-mail: tuganbaev@gmail.com

\textbf{Abstract.} Let $A$ be a ring and $\varphi$ its automorphism. It is proved that skew Laurent series ring $A((x,\varphi ))$ is a right serial ring if and only if $A$ is a right serial right Artinian ring.

The study is supported by Russian Scientific Foundation (project
16-11-10013).

{\bf Kew words.} skew Laurent series ring, right serial ring, Laurent ring

All rings are assumed to be associative and with zero identity element; all modules are unitary and, unless otherwise specified, all modules are right. 

A module is said to be \textit{uniserial} if any two its submodules are comparable with respect to inclusion. Direct sums of uniserial modules are called \textsf{serial} modules. A module $M$ is said to be \textsf{finite-dimensional} if $M$ does not contain an infinite direct sum of non-zero submodules. A module $M$ is said to be \textsf{quotient finite-dimensional} if all factor modules of the module $M$ are finite-dimensional.

Expressions of the form `$A$ is a serial ring' mean that '$A_A$ and $_AA$ are serial modules'.

\textbf{Remark 1.} Let $A$ be a ring and $\varphi$ its automorphism. In \cite[Theorem 2]{Tug08}, it is proved that the skew Laurent series ring $A((x,\varphi ))$ is a right serial ring with the maximum condition on right annihilators if and only if $A$ is a right serial right Artinian ring, if and only if $R$ is a right serial right Artinian ring. 

\textbf{Remark 2.} Let $A$ be a ring and $\varphi$ its automorphism. In \cite[Theorem 1(3)]{Tug16b}, it is proved that skew Laurent series ring $A((x,\varphi ))$ is a serial ring if and only if $A$ is a serial Artinian ring. 

\textbf{Remark 3.} In \cite[Corollary 5.7]{FacH06}, it is proved that any finite direct sum of quotient finite-dimensional modules is a quotient finite-dimensional module; also see \cite[Corollary 5.24]{Fac19}. Therefore, every cyclic right module over a right serial ring is a quotient finite-dimensional module, since any uniserial module is a quotient finite-dimensional module.

\textbf{Remark 4.} In \cite{Tug05}, it is introduced the notion of a Laurent ring $R$ with coefficient ring $A$ and it is proved that the class of all Laurent rings contains all skew Laurent series rings $A((x,\varphi ))$ and pseudo-differential operator rings $A((t^{-1},\delta))$; see \cite[Propositions 7.1, 7.2]{Tug05}. In the case of skew Laurent series rings, the coefficient ring of a Laurent ring coincides with the ordinary coefficient ring $A$ of the ring $R=A((x,\varphi ))$. In addition, in \cite[Proposition 13.5]{Tug05}, it is proved that all cyclic right $R$-modules are finite-dimensional if and only if $R$ is a right Noetherian ring (equivalently, $A$ is a right Noetherian ring).

\textbf{Remark 5.} Let $R$ be a Laurent ring with coefficient ring $A$; see Remark 4. If $R$ is a right serial ring, then it follows from Remarks 3 and 4 that $R$ and $A$ are right Noetherian rings.

The following theorem is the main result of the given paper.

\textbf{Theorem 1.} Let $A$ be a ring and $\varphi$ its automorphism. The following conditions are equivalent. 

\textbf{1)} The skew Laurent series ring $A((x,\varphi ))$ is a right serial ring. 

\textbf{2)} $A((x,\varphi ))$ is a right serial right Artinian ring. 

\textbf{3)} $A$ is a right serial right Artinian ring. 

Theorem 1 follows from Remarks 1 and 5.

The author thanks Alberto Facchini for Remark 3.

\end{document}